\documentclass[]{adamjoucc}
\usepackage{amsmath}
\usepackage{hyperref}
\usepackage{color}
\usepackage{cite}

\begin{document}

\begin{frontmatter}   

\titledata{Reinforcement learning for graph theory, \\ II. Small Ramsey numbers}
{This work was supported and funded by Kuwait University Research Grant No. SM05/22.}

\authordatanoaffil{Mohammad Ghebleh}
{mohammad.ghebleh@ku.edu.kw}
{Corresponding author.},
\authordatanoaffil{Salem Al-Yakoob}
{salem.alyakoob@ku.edu.kw}
{},
\authordatanoaffil{Ali Kanso}
{ali.kanso@ku.edu.kw}
{}\\[1mm]
\affiliation{Department of Mathematics, Faculty of Science, Kuwait University, Safat, Kuwait}

\authordata{Dragan Stevanovic}
{College of Integrative Studies, Abdullah Al Salem University, Khaldiya, Kuwait}
{dragan.stevanovic@aasu.edu.kw}
{On leave from Mathematical Institute of Serbian Academy of Sciences and Arts.}


\begin{abstract}
We describe here how the recent Wagner's approach for applying reinforcement learning to construct examples in graph theory
can be used in the search for critical graphs for small Ramsey numbers.
We illustrate this application by providing lower bounds for the small Ramsey numbers
$R(K_{2,5}, K_{3,5})$, $R(B_3, B_6)$ and $R(B_4, B_5)$ and 
by improving the lower known bound for $R(W_5, W_7)$.
\end{abstract}

\end{frontmatter}

\section{Introduction}

Due to discrete and easily programmable nature of graphs,
computers were used, besides running usual graph algorithms,
to also unearth peculiar examples of graphs answering various research questions
for more than five decades.
Several pieces of software were developed specifically 
to find extremal graphs for arbitrary invariant expressions 
or to conjecture relations between graph invariants.
The most well-known among them are probably 
AutoGraphiX~\cite{agx1,agx2,agx3,agx4} and Graffiti~\cite{graffiti1,graffiti2,graffiti3,graffiti4,graffiti5},
whose Dalmatian heuristic that selectively extends the starting sets of graphs and invariant expressions
was recently reimplemented in~\cite{graffiti6}.

More recently, Adam Wagner~\cite{rl0} has shown how the cross-entropy method, 
a particular reinforcement learning technique developed by Rubinstein~\cite{cema1,cema2},
can be applied to construct graphs that maximize a pre-defined reward,
which he used to find several counterexamples to conjectures in graph theory and combinatorics.
We have reimplemented his approach in~\cite{rl1} in order to make it faster and more stable,
and applied it to disprove a larger set of previously conjectured upper bounds on the Laplacian spectral radius of graphs.

For graphs $G_1,\dots,G_m$,
the $m$-color Ramsey number $R(G_1,\dots,G_m)$ is defined as 
the least integer~$n$ such that
for any edge-coloring of the complete graph~$K_n$ with $m$ colors
there exists some $1\leq i\leq m$
such that the subgraph of~$K_n$ formed by the edges of color~$i$
contains a copy of $G_i$ (not necessarily induced).
An $m$-edge-coloring of~$K_{n}$ is critical
if it does not contain a copy of~$G_i$ having edges of color~$i$ for any $i=1,\dots,m$.
The existence of a critical $m$-edge-coloring of~$K_{n}$
implies that $R(G_1,\dots,G_m)\geq n+1$.
See Radziszowski's dynamic survey~\cite{ramsey} for a wealth of
up-to-date information on the values and bounds for specific Ramsey numbers.

Here we show that, with a small modification,
the cross-entropy method can also be used to find new examples of critical edge-colorings,
using which we improve four lower bounds for small Ramsey numbers.

\section{Cross-entropy method for small Ramsey numbers}

Cross-entropy method is an optimization procedure that employs a feed-forward neural network 
to iteratively construct a new batch of $m$-edge-colorings of~$K_n$ (for fixed $m$ and~$n$)
and to learn the structure of the best performing ones among them.

Edge-colorings are constructed via a sequence of observations and actions.
Each observation is a partial $m$-edge-coloring of~$K_n$ 
together with a one-hot encoding of the current edge that needs to be colored,
while the action is the color selected for the current edge.
The neural network is set to accept an observation as an input
for which, using a few hidden layers, it then outputs the probability distribution of $m$~colors.
The action (i.e., the color of the current edge) is chosen randomly according to this probability distribution,
and the neural network is then used again with an updated observation
to produce the probability distribution of colors for the next edge,
until the ongoing $m$-edge-coloring of~$K_n$ is completed.

This process of constructing $m$-edge-colorings of $K_n$ is repeated for 
a specified number of times,
in order to obtain a full batch of edge-colorings.
Each edge-coloring from the batch receives a reward that is equal to 
the sum of the numbers of copies of the graph~$G_i$ with edges of color~$i$ that it contains for $i=1,\dots,m$.
Obviously, critical edge-colorings are those that have a zero reward.
At this point, a specified percentage of edge-colorings with smallest rewards are selected from the batch,
and the neural network is trained on the sequence of observations and actions that led to these best edge-colorings.
A smaller percentage of the very best edge-colorings is added to the next batch,
to ensure that the best performing edge-colorings are not lost between two batches.

Training of the neural network on the best performing edge-colorings from each batch
ensures that the edge-colorings in the next batch will tend to use actions 
that have previously led to the best perfomers,
which will hopefully decrease the smallest reward in the next batch,
due to the randomized process of selecting edge colors
(the neural network determines the probability distribution of edge colors only).
In any case, keeping a few best survivors from one batch to the next 
ensures that the smallest reward cannot increase.

Over time, the network will start to be retrained on repeated best performers 
and it will eventually almost surely lead to the same edge-colorings over and over again in the later batches,
thus getting stuck at a local minimum.
To postpone this inevitability,
we have introduced a varying percentage of fully randomized actions,
which increases whenever the smallest reward has not decreased for a longer time.
This can help the neural network to avoid learning too specific parts of the structure too soon,
hopefully leading it to find a better local minimum in the long term.
If this local minimum happens to also be a global minimum with a zero reward,
then we have arrived at a critical edge-coloring.
Otherwise, the whole process simply has to be repeated from scratch,
as the weights of the initial neural network are selected entirely randomly,
so that the start of the new training session will very probably lead to a different local minimum.

A Python implementation of this method, with subgraph counting implemented in Java for speedup, 
is available at \url{github.com/dragance106/cema-ramsey}.
It is heavily commented so that the interested user will easily find her way through the code.
Using it, we have been able to find four new critical edge-colorings,
which are presented in the next section.

\section{New lower bounds for Ramsey numbers}

For two graphs $G$ and~$H$,
the join $G+H$ is the graph obtained from the union of $G$ and~$H$
by adding all the edges between a vertex of~$G$ and a vertex of~$H$.
The complete bipartite graph $K_{m,n}$ 
is the join $\overline{K_m}+\overline{K_n}$ of two empty graphs on $m$ and~$n$ vertices, respectively.
The wheel $W_n$ is the join $K_1+C_{n-1}$ of a central vertex and a cycle on $n-1$ vertices.
The book $B_n$ is the join $K_2+\overline{K_n}$, 
which has $n+2$ vertices and can be visualised as $n$~triangular pages having a common edge.

In the following theorems,
we improve the known lower bound $R(W_5,W_7)\geq 13$~\cite{ramsey} and
give the first lower bounds on $R(K_{2,5}, K_{3,5})$, $R(B_3, B_6)$ and $R(B_4, B_5)$ 
by identifying critical edge-colorings using the above cross-entropy method.

\begin{thm}
$R(W_5, W_7)\geq 14$.
\begin{proof}
The following matrix gives a critical edge-coloring of~$K_{13}$ with colors 0 and~1
such that there is 
no $W_5$ with edges of color~0 and no $W_7$ with edges of color~1:
{\small
$$
\begin{bmatrix}       
 &1&1&0&1&0&0&0&1&0&1&1&0\\
1& &0&0&0&1&1&1&0&1&1&1&0\\
1&0& &0&1&1&1&0&0&1&0&0&1\\
0&0&0& &1&1&0&0&1&1&1&1&0\\
1&0&1&1& &0&1&0&0&1&1&0&0\\
0&1&1&1&0& &1&0&0&0&1&0&1\\
0&1&1&0&1&1& &1&1&0&0&0&0\\
0&1&0&0&0&0&1& &1&1&0&0&1\\
1&0&0&1&0&0&1&1& &0&0&1&1\\
0&1&1&1&1&0&0&1&0& &0&1&1\\
1&1&0&1&1&1&0&0&0&0& &0&1\\
1&1&0&1&0&0&0&0&1&1&0& &0\\
0&0&1&0&0&1&0&1&1&1&1&0& 
\end{bmatrix}
$$
}

Alternatively, the graph determined by this adjacency matrix (edge color~1) does not contain $W_7$ as a subgraph,
while its complement (edge color~0) does not contain $W_5$ as a subgraph.
\end{proof}
\end{thm}

\begin{thm}
$R(K_{2,5}, K_{3,5})\geq 20$.
\begin{proof}
The following matrix gives a critical edge-coloring of~$K_{19}$ with colors 0 and~1
such that there is 
no $K_{2,5}$ with edges of color~0 and no $K_{3,5}$ with edges of color~1:
{\small
$$
\begin{bmatrix}
 &0&1&1&1&1&0&0&1&1&1&0&1&1&0&0&1&0&0\\
0& &1&0&0&1&1&1&0&1&1&1&0&1&1&1&0&0&1\\
1&1& &0&0&1&1&1&1&1&0&1&0&1&0&0&0&1&0\\
1&0&0& &1&1&0&0&1&1&1&1&1&0&1&1&0&1&0\\
1&0&0&1& &0&1&0&0&1&0&1&1&1&1&0&0&1&1\\
1&1&1&1&0& &0&1&1&0&1&1&0&0&1&0&0&0&1\\
0&1&1&0&1&0& &0&1&1&0&1&1&1&0&1&0&0&1\\
0&1&1&0&0&1&0& &1&0&0&0&1&1&1&0&1&1&1\\
1&0&1&1&0&1&1&1& &0&0&0&1&0&0&1&1&0&1\\
1&1&1&1&1&0&1&0&0& &0&0&0&1&0&1&1&1&1\\
1&1&0&1&0&1&0&0&0&0& &1&1&1&1&1&1&0&0\\
0&1&1&1&1&1&1&0&0&0&1& &0&1&1&0&0&1&0\\
1&0&0&1&1&0&1&1&1&0&1&0& &1&1&1&0&0&0\\
1&1&1&0&1&0&1&1&0&1&1&1&1& &0&0&1&0&0\\
0&1&0&1&1&1&0&1&0&0&1&1&1&0& &1&0&1&1\\
0&1&0&1&0&0&1&0&1&1&1&0&1&0&1& &1&1&1\\
1&0&0&0&0&0&0&1&1&1&1&0&0&1&0&1& &1&1\\
0&0&1&1&1&0&0&1&0&1&0&1&0&0&1&1&1& &0\\
0&1&0&0&1&1&1&1&1&1&0&0&0&0&1&1&1&0& 
\end{bmatrix}
$$
}
\end{proof}       
\end{thm}

\begin{thm}       
$R(B_3, B_6)\geq 17$.
\begin{proof}
The following matrix gives a critical edge-coloring of~$K_{16}$ with colors 0 and~1
such that there is 
no $B_3$ with edges of color~0 and no $B_6$ with edges of color~1:
{\small
$$
\begin{bmatrix}       
 &1&1&1&1&0&1&1&0&0&1&1&0&1&0&0\\
1& &0&0&1&1&1&1&1&1&0&1&1&0&0&0\\
1&0& &1&1&0&1&0&0&0&0&1&1&1&0&1\\
1&0&1& &0&0&0&1&1&0&1&1&1&1&1&0\\
1&1&1&0& &0&1&0&1&1&1&0&1&0&0&1\\
0&1&0&0&0& &1&1&1&1&0&1&0&1&1&1\\
1&1&1&0&1&1& &1&1&0&1&0&0&0&1&0\\
1&1&0&1&0&1&1& &0&1&1&0&1&0&1&1\\
0&1&0&1&1&1&1&0& &0&1&1&1&0&1&1\\
0&1&0&0&1&1&0&1&0& &1&1&1&1&0&1\\
1&0&0&1&1&0&1&1&1&1& &0&0&1&1&0\\
1&1&1&1&0&1&0&0&1&1&0& &1&1&0&0\\
0&1&1&1&1&0&0&1&1&1&0&1& &0&0&1\\
1&0&1&1&0&1&0&0&0&1&1&1&0& &1&1\\
0&0&0&1&0&1&1&1&1&0&1&0&0&1& &1\\
0&0&1&0&1&1&0&1&1&1&0&0&1&1&1& 
\end{bmatrix}
$$
}
\end{proof}
\end{thm}

\begin{thm}       
$R(B_4, B_5)\geq 18$.
\begin{proof}
The following matrix gives a critical edge-coloring of~$K_{17}$ with colors 0 and~1
such that there is 
no $B_4$ with edges of color~0 and no $B_5$ with edges of color~1:
{\small
$$
\begin{bmatrix}       
 &0&0&0&0&1&1&1&0&0&1&1&0&1&1&1&0\\
0& &0&1&0&1&0&0&1&1&1&1&1&1&0&1&0\\
0&0& &0&1&0&0&1&1&1&1&0&0&1&0&1&1\\
0&1&0& &1&0&1&0&0&0&0&1&1&1&0&1&1\\
0&0&1&1& &1&1&0&1&0&1&0&0&1&0&0&1\\
1&1&0&0&1& &0&0&1&0&1&0&1&0&1&1&1\\
1&0&0&1&1&0& &0&0&1&1&1&0&1&1&0&1\\
1&0&1&0&0&0&0& &1&1&0&1&1&0&1&1&0\\
0&1&1&0&1&1&0&1& &0&1&0&1&1&0&0&0\\
0&1&1&0&0&0&1&1&0& &1&1&1&0&1&0&1\\
1&1&1&0&1&1&1&0&1&1& &1&0&0&0&0&0\\
1&1&0&1&0&0&1&1&0&1&1& &0&0&0&1&0\\
0&1&0&1&0&1&0&1&1&1&0&0& &0&1&0&1\\
1&1&1&1&1&0&1&0&1&0&0&0&0& &1&1&0\\
1&0&0&0&0&1&1&1&0&1&0&0&1&1& &0&1\\
1&1&1&1&0&1&0&1&0&0&0&1&0&1&0& &1\\
0&0&1&1&1&1&1&0&0&1&0&0&1&0&1&1& 
\end{bmatrix}
$$
}
\end{proof}
\end{thm}

\section{Concluding remarks}

Numerous experiments with the implemented method showed that 
higher learning percentages and larger batch sizes tend to help it
to more easily identify promising directions for learning.
While the results from the previous section do prove that 
the cross-entropy method may lead to new examples of critical edge-colorings for small Ramsey numbers,
most often it actually leaves the questions simply unanswered.

As a prominent example,
the best 2-edge-colorings of~$K_{17}$ that it managed to produce after hundreds of runs, 
have between 1 and~8 monochromatic copies of $B_3$ and~$B_6$, depending on the run.
Regardless of the number of runs,
one still cannot be certain whether the reason for this is because truly $R(B_3,B_6)=17$
or we just did not attempt that one additional lucky run of the cross-entropy method?

Another question is how to predict when is the time to end the current run?
We have observed cases when the reward reaches a plateau after 300--400 batches and 
then it remains constant until a sudden decrease more than 10,000 batches later.
This decrease was apparently a consequence of a lucky randomization of actions that accidentally led to a smaller reward,
moving away from the local minimum where it was stuck for thousands of batches.
However, what this learning method really needs is an improved way of getting out of local minima.
At present, a jump to a better local minimum through randomization of actions has to be made within a single batch,
as otherwise such randomized colorings will not qualify within the surviving or even the learning percentage,
and they will be already discarded by the time the next batch is constructed.

The final remark is that, 
although the cross-entropy method does manage to construct some critical edge-colorings,
it may not actually be the best way for doing it.
Namely, each critical edge-coloring of~$K_n$ automatically yields a critical edge-coloring of~$K_{n-1}$
simply by deleting an arbitrary vertex of~$K_n$.
Thus, it may actually be more beneficial to proceed inductively with learning:
assuming the machine learning model~$M_{n-1}$ knows how to construct critical edge-colorings of~$K_{n-1}$
to then somehow use the acquired knowledge of~$M_{n-1}$
to train the model~$M_n$ that will also construct critical edge-colorings of~$K_n$,
provided $n<R(G_1,\dots,G_m)$.

We leave the resolution of these questions for future work.

\end{document}